\documentclass[12pt,a4paper,twoside]{article}

\newcommand{\pageformat}[6]{\setlength{\hoffset}{-1in}
                  \setlength{\voffset}{-1in}
                  \addtolength{\hoffset}{#5}
                            \addtolength{\voffset}{#6}
                            \setlength{\oddsidemargin}{#1}
                            \setlength{\evensidemargin}{#2}
                            \setlength{\textwidth}{\paperwidth}
                  \addtolength{\textwidth}{-\oddsidemargin}
                  \addtolength{\textwidth}{-\evensidemargin}
                  \addtolength{\textwidth}{-\marginparsep}
                  \addtolength{\textwidth}{-\marginparwidth}
                            \setlength{\topmargin}{#3}
                            \setlength{\textheight}{\paperheight}
                  \addtolength{\textheight}{-\topmargin}
                  \addtolength{\textheight}{-\headheight}
                  \addtolength{\textheight}{-\headsep}
                  \addtolength{\textheight}{-\footskip}
                  \addtolength{\textheight}{-#4}}
\pageformat{2cm}{3cm}{24mm}{24mm}{1pt}{0pt}

\usepackage{ifthen}
\newboolean{article}
    \setboolean{article}{true}
\newboolean{report}
\newboolean{book}
\newboolean{letter}
\newboolean{german}
\newboolean{italian}
\newboolean{nobaselinestretch}
\newboolean{nosectionappendix}
\newboolean{oldtoc}
\newboolean{nosectionequn}
\newboolean{notheorem}

\ifthenelse{\boolean{german}}{
    \usepackage{german}}{}

\usepackage[latin1]{inputenc}

\usepackage{amsmath}
\usepackage{amssymb}
\usepackage[mathscr]{eucal}

\ifthenelse{\boolean{notheorem}}{}{
    \usepackage{theorem}}

\ifthenelse{\boolean{nobaselinestretch}}{}{
    \renewcommand{\baselinestretch}{1.25}}

\newenvironment{env}[2]{\begin{#1}#2\end{#1}}{}
    \newcommand{\beq}[1]{\begin{env}{equation}{#1}}
    \newcommand{\beqn}[1]{\begin{env}{equation*}{#1}}
    \newcommand{\bal}[1]{\begin{env}{align}{#1}}
    \newcommand{\baln}[1]{\begin{env}{align*}{#1}}
    \newcommand{\bga}[1]{\begin{env}{gather}{#1}}
    \newcommand{\bgan}[1]{\begin{env}{gather*}{#1}}
    \newcommand{\bflal}[1]{\begin{env}{flalign}{#1}}
    \newcommand{\bflaln}[1]{\begin{env}{flalign*}{#1}}
    \newcommand{\bmu}[1]{\begin{env}{multline}{#1}}
    \newcommand{\bmun}[1]{\begin{env}{multline*}{#1}}
    \newcommand{\bsp}[1]{\begin{env}{split}{#1}}

    \newcommand{\eeq}{\end{env}}
    \newcommand{\eeqn}{\end{env}}
    \newcommand{\eal}{\end{env}}
    \newcommand{\ealn}{\end{env}}
    \newcommand{\ega}{\end{env}}
    \newcommand{\egan}{\end{env}}
    \newcommand{\eflal}{\end{env}}
    \newcommand{\eflaln}{\end{env}}
    \newcommand{\emu}{\end{env}}
    \newcommand{\emun}{\end{env}}
    \newcommand{\esp}{\end{env}}

\newcommand{\lf}{\vspace{2ex}}

\renewcommand{\bf}[1]{\textbf{#1}}
\renewcommand{\it}[1]{\textit{#1}}

\renewcommand{\sf}[1]{\textsf{#1}}

\renewcommand{\tt}[1]{\texttt{#1}}
\newcommand{\hl}[1]{\bf{\it{#1}}}

\newcommand{\mbf}[1]{\mathbf{#1}}
\newcommand{\msf}[1]{\text{\small$\sf{#1}$}}

\newcommand{\cmc}[1]{\mathcal{#1}}
\newcommand{\eus}[1]{\mathscr{#1}}
\newcommand{\euf}[1]{\mathfrak{#1}}
\newcommand{\bb}[1]{\mathbb{#1}}

\newcommand{\nbd}[1]{$#1$\nobreakdash--}
\newcommand{\ol}[1]{\overline{#1}}

\newcommand{\wt}[1]{\widetilde{#1}}

\newcommand{\vt}{\vartheta}

\newcommand{\om}{\omega}

\newcommand{\bfam}[1]{\bigl(#1\bigr)}

\newcommand{\AB}[1]{\langle#1\rangle}
\newcommand{\bAB}[1]{\bigl\langle#1\bigr\rangle}

\newcommand{\CB}[1]{\{#1\}}
\newcommand{\bCB}[1]{\bigl\{#1\bigr\}}
\newcommand{\BCB}[1]{\Bigl\{#1\Bigr\}}

\newcommand{\RO}[1]{[#1)}

\newcommand{\set}[2][]{
    \ifthenelse{\equal{#1}{}}{
        \CB{#2}}{
        \CB{#1~|~#2}}}
\newcommand{\bset}[2][]{
    \ifthenelse{\equal{#1}{}}{
        \bCB{#2}}{
        \bCB{#1~|~#2}}}
\newcommand{\Bset}[2][]{
    \ifthenelse{\equal{#1}{}}{
        \BCB{#2}}{
        \BCB{#1~\big|~#2}}}

\DeclareMathOperator{\ls}{\normalfont\msf{span}}
\DeclareMathOperator{\cls}{\ol{\ls}}
\DeclareMathOperator*{\limind}{lim\,ind}

\DeclareMathOperator{\id}{\normalfont\msf{id}}

\newcommand{\N}{\bb{N}}

\newcommand{\R}{\bb{R}}

\newcommand{\T}{\bb{T}}

\newcommand{\cB}{\cmc{B}}

\newcommand{\sB}{\eus{B}}

\newcommand{\sK}{\eus{K}}

\newcommand{\eH}{\euf{H}}

\newcommand{\U}{\mbf{1}}

\newcommand{\G}{\Gamma}

\ifthenelse{\boolean{nosectionequn}}{}{
    \numberwithin{equation}{section}
    }

\ifthenelse{\boolean{article}\or\boolean{letter}\or\boolean{nosectionequn}}{
    \setboolean{nosectionappendix}{true}}{}
\ifthenelse{\boolean{nosectionappendix}}{}{
    \renewcommand{\appendix}{
        \chapter*{\appendixname}
        \addcontentsline{toc}{chapter}{\appendixname}
        \renewcommand{\thesection}{\Alph{section}}
        \setcounter{section}{0}}}
   
\ifthenelse{\boolean{report}\or\boolean{book}}{
    }{}

\ifthenelse{\boolean{notheorem}}{}{
        \newcommand{\mnname}{Mathematical note.}
        \newcommand{\enname}{End of the note.}
        \newcommand{\definame}{Definition.}
        \newcommand{\propname}{Proposition.}
        \newcommand{\lemname}{Lemma.}
        \newcommand{\exname}{Example.}
        \newcommand{\exername}{Exercise.}
        \newcommand{\remname}{Remark.}
        \newcommand{\obname}{Observation.}
        \newcommand{\thmname}{Theorem.}
        \newcommand{\corname}{Corollary.}
        \newcommand{\proofname}{Proof.}
    \ifthenelse{\boolean{german}}{
        \renewcommand{\mnname}{Mathematische Notiz.}
        \renewcommand{\enname}{Ende der Notiz.}
        \renewcommand{\exname}{Beispiel.}
        \renewcommand{\exername}{Übung.}
        \renewcommand{\remname}{Bemerkung.}
        \renewcommand{\obname}{Beobachtung.}
        \renewcommand{\thmname}{Satz.}
        \renewcommand{\corname}{Korollar.}
        \renewcommand{\proofname}{Beweis.}}{}
    \ifthenelse{\boolean{italian}}{
        \renewcommand{\mnname}{Nota matematica.}
        \renewcommand{\enname}{Fina della nota.}
        \renewcommand{\definame}{Definizione.}
        \renewcommand{\propname}{Proposizione.}
        \renewcommand{\exname}{Esempio.}
        \renewcommand{\exername}{Esercizio.}
        \renewcommand{\remname}{Nota.}
        \renewcommand{\obname}{Osservazione.}
        \renewcommand{\thmname}{Teorema.}
        \renewcommand{\corname}{Corollario.}
        \renewcommand{\proofname}{Dimostrazione.}

       \renewcommand{\appendixname}{Appendice}

       }{}
    \theoremheaderfont{\normalfont\bfseries}
    \theoremstyle{change}
        \theorembodyfont{\rmfamily}
            \newtheorem{emp}{}[section]
                \newcommand{\bemp}[1][]{
                    \begin{emp}\hskip-\labelsep\bf{#1}\hskip\labelsep}
                \newcommand{\eemp}{\end{emp}}
\newtheorem{itemp}[emp]{}
                \newcommand{\bitemp}[1][]{
                    \begin{itemp}\hskip-\labelsep\bf{#1}\hskip\labelsep\normalfont\itshape}
                \newcommand{\eitemp}{\end{itemp}}
            \newtheorem{mn}[emp]{\mnname}
                \newcommand{\bnm}{\begin{mn}~\begin{quotation}\renewcommand{\baselinestretch}{1}\small\noindent\ignorespaces}
                \newcommand{\enm}{\end{quotation}\hfill\bf{\enname}\end{mn}}
            \newtheorem{ex}[emp]{\exname}
                \newcommand{\bex}{\begin{ex}}
                \newcommand{\eex}{\end{ex}}
            \newtheorem{exer}[emp]{\exername}
                \newcommand{\bexer}{\begin{exer}}
                \newcommand{\eexer}{\end{exer}}
            \newtheorem{defi}[emp]{\definame}
                \newcommand{\bdefi}{\begin{defi}}
                \newcommand{\edefi}{\end{defi}}
            \newtheorem{rem}[emp]{\remname}
                \newcommand{\brem}{\begin{rem}}
                \newcommand{\erem}{\end{rem}}
            \newtheorem{ob}[emp]{\obname}
                \newcommand{\bob}{\begin{ob}}
                \newcommand{\eob}{\end{ob}}
        \theorembodyfont{\normalfont\itshape}
            \newtheorem{thm}[emp]{\thmname}
                \newcommand{\bthm}{\begin{thm}}
                \newcommand{\ethm}{\end{thm}}
            \newtheorem{prop}[emp]{\propname}
                \newcommand{\bprop}{\begin{prop}}
                \newcommand{\eprop}{\end{prop}}
            \newtheorem{cor}[emp]{\corname}
                \newcommand{\bcor}{\begin{cor}}
                \newcommand{\ecor}{\end{cor}}
            \newtheorem{lem}[emp]{\lemname}
                \newcommand{\blem}{\begin{lem}}
                \newcommand{\elem}{\end{lem}}
\newenvironment{empn}[1]{\lf\noindent\bf{#1}\ignorespaces\hskip\labelsep}{\lf}
		\newcommand{\bempn}[1]{\begin{empn}{#1}}
		\newcommand{\eempn}{\end{empn}}
		\newcommand{\bitempn}[1]{\begin{empn}{#1}\normalfont\itshape}
		\newcommand{\eitempn}{\end{empn}}
                \newcommand{\bnmn}{\begin{empn}{\mnname}~\begin{quotation}\renewcommand{\baselinestretch}{1}\small\noindent\ignorespaces}
                \newcommand{\enmn}{\end{quotation}\hfill\bf{\enname}\end{empn}}
		\newcommand{\bexn}{\begin{empn}{\exname}}
		\newcommand{\eexn}{\end{empn}}
		\newcommand{\bexern}{\begin{empn}{\exername}}
		\newcommand{\eexern}{\end{empn}}
		\newcommand{\bdefin}{\begin{empn}{\definame}}
		\newcommand{\edefin}{\end{empn}}
		\newcommand{\bremn}{\begin{empn}{\remname}}
		\newcommand{\eremn}{\end{empn}}
		\newcommand{\bobn}{\begin{empn}{\obname}}
		\newcommand{\eobn}{\end{empn}}

\newcommand{\qedsymbol}{~\rule[-0.35mm]{2mm}{2mm}}
    \newcounter{proof}[emp]
    \newenvironment{Proof}[1]{
        \vspace{1ex}
        \renewcommand{\item}[1][\stepcounter{proof}(\roman{proof})]%
            {##1\hskip\labelsep}
        \noindent\textsc{#1\hskip\labelsep}}{
        \nolinebreak\qedsymbol}
    \newcommand{\proof}[1][\proofname]{
        \begin{Proof}{#1}\ignorespaces}
    \newcommand{\qed}{\end{Proof}}
    \newcommand{\noqed}{
        \renewcommand{\qedsymbol}{}
        \end{Proof}}}
    \ifthenelse{\boolean{italian}}{
        \renewcommand{\proofname}{Dimostrazione.}}{}

\usepackage{txfonts}

\usepackage[matrix,arrow,curve]{xy}

\begin{document}

\title{Spatial $E_0$--Semigroups are Restrictions\\of Inner Automorphismgroups}
\author{}
\author{
~\\
Michael Skeide\thanks{This work is supported by research fonds of the Department S.E.G.e S.\ of University of Molise.}
}
\date{}

\maketitle

\noindent
If $\vt$ is a strict \nbd{E_0}semigroup on some $\sB^a(E)$, then Skeide \cite{Ske02,Ske04p} and Muhly, Skeide and Solel \cite{MSS03p} associate with $\vt$ a product system of correspondences. We say the \nbd{E_0}semigroup is \hl{spatial}, if the associated product system is spatial in the sense of Skeide \cite{Ske01p}. The main goal of these notes is to establish the following theorem that just restates the title of these notes in a more specific form. Further terminology used in the theorem will be discussed after the notes and the example.

\lf\noindent
\bf{Main theorem.}
\it{
Suppose that $E_+$ is a Hilbert module over a (unital) \nbd{C^*}algebra $\cB$ and that $\vt=\bfam{\vt_t}_{t\in\T}$ is a spatial strict \nbd{E_0}semigroup on $\sB^a(E_+)$. Then there exists a correspondence $E_-$ over $\cB$ and a semigroup $w=\bfam{w_t}_{t\in\T}$ of unitaries $w_t$ on $E:=E_+\odot E_-$ such that the canonical homomorphism $\sB^a(E_+)\rightarrow\sB^a(E_+)\odot\id_{E_-}$ is an isomorphism and such that for every $t\in\T$ the restriction of $\alpha_t:=w_t\bullet w_t^*$ to $\sB^a(E_+)\odot\id_{E_-}$ is $\vt_t\odot\id_{E_-}$.
}

\lf\noindent
\bf{Notes.}
1.)
The result has an obvious variation for normal \nbd{E_0}semigroups when $E_+$ is a von Neumann (or \nbd{W^*}) module. Just replace correspondences by von Neumann (or \nbd{W^*}) correspondences and their tensor products.

2.)
Following the lines of Skeide \cite{Ske03b} it is easy to show that the unitary semigroup $w$ reflects the continuity properties of the \nbd{E_0}semigroup $\vt$. But we do not have enough space to include the proof of such technicalities here.

3.)
Unfortunately, for Hilbert modules the condition that $\sB^a(E_+)\rightarrow\sB^a(E_+)\odot\id_{E_-}$ injective is not automatic, if the the left action of $\cB$ on $E_-$ fails to be faithful.

\lf\noindent
\bf{Example.~~}
The usual time-shift endomorphism (CCR-flow) on the symmetric Fock space $\G(L^2(\R_+,K))$ may be understood as the restriction of of the unitarily implemented time-shift automorphism on $\G(L^2(\R,K))=\G(L^2(\R_+,K))\otimes\G(L^2(\R_-,K))$ to $\sB\bfam{\G(L^2(\R_+,K))}\otimes\id_{\G(L^2(\R_-,K))}$. The same is true for \it{time-ordered} Fock modules \cite{BhSk00}, the module analogue of the symmetric Fock space.

\lf
Now we explain in detail the terms used in the theorem. The semigroup with identity $\T$ is either $\R_+=\RO{0,\infty}$ or $\N_0=\CB{0,1,\ldots}$. A \hl{correspondence} over $\cB$ (or a Hilbert \nbd{\cB}bimodule) is a Hilbert \nbd{\cB}module with a nondegenerate left action of $\cB$ as a representation by adjointable operators. $\sB^a(E)$ denotes the \nbd{C^*}algebra of all adjointable operators on a Hilbert module $E$. The tensor product is the internal tensor product over $\cB$. An \hl{\nbd{E_0}semigroup} is a semigroup of unital endomorphisms, and a unital endomorphism of $\sB^a(E_+)$ is \hl{strict} if its restriction to the \hl{compacts} $\sK(E_+):=\cls\CB{xy^*~(x,y\in E_+)}$ acts nondegenerately on $E_+$, where $xy^*$ denotes the \hl{rank-one operator} $z\mapsto x\AB{y,z}$. (This is the simplest and most useful criterion for strictness; see, for instance, Lance \cite{Lan95}.)

A \hl{product system} is a family $E^\odot=\bfam{E_t}_{t\in\T}$ of correspondences over $\cB$ with $E_0=\cB$ and an associative family of isomorphisms (that is, bilinear unitaries) $u_{st}\colon E_s\odot E_t\rightarrow E_{s+t}$ (being the canonical ones when $s=0$ or $t=0$). Using the representation theory of \cite{MSS03p}, \cite{Ske04p} associates with every $\vt_t$ $(t>0)$ the correspondence $E_+^*\odot{~}_{\vt_t}E_+$, where $E_+^*$ is the \hl{dual} correspondence from $\cB$ to $\sB^a(E_+)$ of $E_+$ with inner product $\AB{x^*,y^*}=xy^*$ and obvious \nbd{\cB}\nbd{\sB^a(E_+)}bimodule operations, while $_{\vt_t}E_+$ is $E_+$ viewed as correspondence from $\sB^a(E_+)$ to $\cB$ with left action of $\sB^a(E_+)$ via $\vt_t$. (If $E_+$ is \hl{full}, that is, if $\cls\AB{E_+,E_+}=\cB$, then this definition applies also for $t=0$. Otherwise, we have to put $E_0=\cB$ by hand.) The isomorphisms $u_{st}$ are determined by $u_{st}((x_s^*\odot_s x'_s)\odot(y_t^*\odot_t y'_t))=x_s^*\odot_{s+t}\vt_t(x'_sy_t^*)y'_t$, where by $\odot_t$ we indicate that the tensor product is that of $E_+^*\odot{~}_{\vt_t}E_+=E_t$. Product system and \nbd{E_0}semigroup are related by the family $u=\bfam{u_t}_{t\in\T}$ of unitaries $u_t\colon E\odot E_t\rightarrow E$ determined by $u_t(x\odot(y^*\odot_t z)=\vt_t(xy^*)z$ which give back $\vt_t$ as $\vt_t(a)=u_t(a\odot\id_{E_t})u_t^*$.

A \hl{unit} for a product system $E^\odot=\bfam{E_t}_{t\in\T}$ is a family $\xi^\odot=\bfam{\xi_t}_{t\in\T}$ of vectors $\xi_t\in E_t$ such that $u_{st}(\xi_s\odot\xi_t)=\xi_{s+t}$ and $\xi_0=\U\in E_0=\cB$. This implies necessarily that $\cB$ is unital. (As observed in Bhat and Skeide \cite{BhSk00} the definition of the inner product in internal tensor products implies that the mappings $T_t:=\AB{\xi_t,\bullet\xi_t}$ form a CP-semigroup on $\cB$. Without the condition on $\xi_0$, the mapping $T_0$ could never be the identity. See also the discussion of nonunital $\cB$ in Skeide \cite{Ske04p}.) According to \cite{Ske01p} a product system is \hl{spatial} if it admits a central unital unit $\om^\odot=\bfam{\om_t}_{t\in\T}$. Here \hl{central} means that $b\om_t=\om_tb$ for all $t\in\T,b\in\cB$, and \hl{unital} means that $\AB{\om_t,\om_t}=\U$ for all $t\in\T$. Spatiality of the product system of $\vt$ implies that $E_+$ is full. Of course, a central unital unit $\om^\odot$ generates the \hl{trivial} CP-semigroup $\AB{\om_t,\bullet\om_t}=\id_\cB$.

\lf\noindent
\bf{Remark.}
For $\cB$ unital and $E_+$ full one may show that spatiality of the \nbd{E_0}semigroup $\vt$ is equivalent to existence of a semigroup $v=\bfam{v_t}_{t\in\T}$ of intertwining isometries $v_t\in\sB^a(E_+)$ for $\vt$, that is, $\vt_t(a)v_t=v_ta$. (This is closest to Powers' original definition of spatial \nbd{E_0}semigroups \cite{Pow87} in the case when $E_+$ is just a Hilbert space.) In fact, if there is a central unital unit $\om^\odot$, then $v_t:=u_t(\id_{E_+}\odot\om_t)$ defines such a semigroup of intertwining isometries. (Here $\id_{E_+}\odot\om_t$ denotes the mapping $x\mapsto x\odot\om_t$.) Conversely, given such a semigroup $v$ by \it{general abstract nonsense} one may show the converse. (This involves Rieffel's fundamental results on \it{Morita equivalence} \cite{Rie74,Rie74a} together with such simple observations like $E_+^*$ (full!) is a Morita equivalence from $\cB$ to $\sK(E_+)$ and $\sB^a(\cB,E_+)=E_+$ ($\cB$ unital!) thinking of $x\in E_+$ as the mapping $b\mapsto xb$.)

\lf
We will also show a supplement to the main theorem regarding weak dilations in the sense of \cite{BhSk00}. The pair $(E_+,\vt)$ is a \hl{weak dilation} of a (necessarily unital) CP-semigroup, if there exists a unit vector $\xi_+\in E_+$ (\hl{unit vector} means that $\AB{\xi,\xi}=\U$) such that the projection $p_0:=\xi_+\xi_+^*$ is \hl{increasing} for $\vt$, that is, $\vt_t(p_0)\ge p_0$ for all $t\in\T$. In this case $T_t(b):=\AB{\xi_+,\vt_t(\xi_+b\xi_+^*)\xi_+}$, indeed, defines a unital CP-semigroup $T$ and $\vt$ , clearly, is a dilation of $T$ (under the embedding $b\mapsto\xi_+b\xi_+^*$). Also, it is not difficult to check that under these circumstances $E^\odot$ has a unital unit $\xi^\odot$ such that $T_t=\AB{\xi_t,\bullet\xi_t}$; see \cite{Ske02}.

\lf\noindent
\bf{Supplement.}
\it{
If $(E_+,\vt,\xi_+)$ is a weak dilation of the unital CP-semigroup $T$, then the correspondence $E_-$ in the main theorem can be chosen such that $E_-$ contains a (central) unit vector $\om_-$ and the semigroup $\alpha$ is a dilation of $T$ with respect to the embedding $b\mapsto\xi b\xi^*\odot\id_{E_-}$ in the vector expectation $\AB{\xi_+\odot\om_-,\bullet\xi_+\odot\om_-}$.
}

\lf
We note that $(E_+\odot E_-,\alpha,\xi_+\odot\om_-)$ is a weak dilation of $T$, if and only if $T$ is the trivial CP-semigroup and if the projection $\id_{E_+}\odot\om_-\om_-^*$ is invariant for $\alpha$. It \hl{is}, generally, a weak dilation with respect to the filtration $\alpha_t(p_0\odot\id_{E_-})=\vt_t(p_0)\odot\id_{E_-}$, if we apply the weaker hypothesis of Bhat and Parthasarathy \cite{BhPa94} (rephrased suitably in terms of Hilbert modules).

\section{Units and inductive limits}\label{ILsec}

As a motivation for the construction of $E_-$ we repeat an inductive limit construction from \cite{BhSk00} that reverses in some sense the construction of the product system $E^\odot $ from $\vt$ in the case when $(E_+,\vt)$ is weak dilation.

So let $E^\odot$ be a product system and $\xi^\odot$ a unital unit for $E^\odot$. Then $\xi_s\odot\id_{E_t}\colon x_t\mapsto u_{st}(\xi_s\odot x_t)$ defines an isometry in $\sB^a(E_t,E_{s+t})$. The family $\bfam{E_t}_{t\in\T}$ together with the family of embeddings $\xi_s\odot\id_{E_t}$ forms an inductive system of Hilbert \nbd{\cB}modules (not of correspondences!) and the completion of the algebraic inductive limit is Hilbert \nbd{\cB}module which we denote by $E_\infty$.

The factorization $E_s\odot E_t=E_{s+t}$ (we surpress the mappings $u_{st}$) survives the inductive limit and gives rise to a factorization $E_\infty\odot E_t=E_\infty$ fulfilling the associativity condition $(E_\infty\odot E_s)\odot E_t=E_\infty\odot(E_s\odot E_t)$. Therefore, $\wt{\vt}_t(a):=a\odot\id_{E_t}$ defines a strict \nbd{E_0}semigroup on $\sB^a(E_\infty)$.

As $\xi_s\odot\xi_t=\xi_{s+t}$, the unit vector $\xi_t\in E_t$, when embedded into $E_{s+t}$, coincides with $\xi_{s+t}$. Therefore, $E_\infty$ contains a distinguished unit vector $\xi_\infty$, the inductive limit of all the $\xi_t$, and $(E_\infty,\wt{\vt},\xi_\infty)$ is a weak dilation of the CP-semigroup $T_t=\AB{\xi_t,\bullet\xi_t}$.

\brem
If $(E_+,\vt,\xi_+)$ is a weak dilation and if $E^\odot$ and $\xi^\odot$ are product system and unit associated with that dilation, then $E_\infty$ is identified naturally as the submodule $\ol{\bigcup}_{t\in\T}\vt_t(p_0)E_+$ of $E_+$. The dilation is \hl{primary}, if $E_+=E_\infty$. In this case also $\vt=\wt{\vt}$ and $\xi_+=\xi_\infty$.
\erem

\section{Proof of the main theorem and its supplement}

Already in Skeide \cite{Ske01} we noted that in the case of a central unital unit $\om^\odot$ the preceding inductive limit can be performed into \it{the other direction}, that is, using embeddings $\id_{E_t}\odot\om_s$ rather than $\om_s\odot\id_{E_t}$. Indeed, in the identification $E_t=\sB^a(\cB,E_t)$ ($\cB$ is unital!), the mapping $\om_t\colon b\mapsto\om_tb$ is actually bilinear and, therefore, can be amplificated as a right factor in a tensor product.

The embeddings $E_t\rightarrow E_{s+t}$ we obtain in that way are, indeed, bilinear so that we obtain an inductive limit $E_-$ which is a bimodule, that is, a correpondence over $\cB$. We have now a family of bilinear unitaries $u_t^-\colon E_t\odot E_-\rightarrow E_-$. Clearly, $\om_-:=\limind_{t\in\T}\om_t$ is a central unit vector, that is in particular, $\AB{\om,b\om}=b$ for all $b\in\cB$.

Now we can put together $u_t$ and $u_t^-$ to form a unitary $w_t$ on $E_+\odot E_-$. We define
\beqn{
w_t=(u_t\odot\id_{E_-})(\id_{E_+}\odot{u_t^-}^*).
}\eeqn
What $w_t$ does is simply
\beqn{
E_+\odot E_-
~\longrightarrow~
E_+\odot(E_t\odot E_-)
~\longrightarrow~
(E_+\odot E_t)\odot E_-
~\longrightarrow~
E_+\odot E_-
}\eeqn
or
\beqn{
w_t(x\odot u_t^-(y_t\odot z))
~=~
u_t(x\odot y_t)\odot z.
}\eeqn
Thinking of $u_t$ as identification $E_+\odot E_t=E_+$, as we did in many papers, and of $u_t^-$ as identification $E_t\odot E_-=E_-$, then $w_t$ is just rebracketting
\beqn{
x\odot\underbrace{y_t\odot z}_{E_-}
~\longmapsto~
\underbrace{x\odot y_t}_{E_+}\odot z.
}\eeqn
Clearly, the $w_t$ define a semigroup. And for every $a\in\sB^a(E_+)$, as $(a\odot\id_{E_-})(\id_{E_+}\odot u_t^-)=(\id_{E_+}\odot u_t^-)(a\odot\id_{E_t}\odot\id_{E_-})$, we find
\bmun{
\alpha_t(a\odot\id_{E_-})
~=~
(u_t\odot\id_{E_-})(\id_{E_+}\odot{u_t^-}^*)(a\odot\id_{E_-})(\id_{E_+}\odot{u_t^-})(u_t^*\odot\id_{E_-})
\\
~=~
(u_t\odot\id_{E_-})(a\odot\id_{E_t}\odot\id_{E_-})(u_t^*\odot\id_{E_-})
~=~\vt_t(a)\odot\id_{E_-}.
}\emun
This proves the main theorem.

The assertions of the supplement follow simply from the preceding calculation by observing that $\bAB{(\xi_+\odot\om_-),(\vt_t(a)\odot\id_{E_-})(\xi_+\odot\om_-)}=\bAB{\om_-,\AB{\xi_+,\vt_t(a)\xi_+}\om_-}=\AB{\xi_+,\vt_t(a)\xi_+}$.

\brem
For proving the note after the supplement, it is sufficient to observe, that $(\xi_+\odot\om_-)(\xi_+\odot\om_-)^*$ is increasing for $\alpha$, if and only if $\AB{\om_t,\xi_t}=\U$ for all $t\in\T$. This implies $\xi_t\xi_t^*\ge\om_t\om_t^*$ and $\xi_t\xi_t^*\le\om_t\om_t^*$, hence, $\xi_t\xi_t^*=\om_t\om_t^*$ and from this $\xi_t=\xi_t\xi_t^*\xi_t=\om_t\om_t^*\xi_t=\om_t\U=\om_t$, that is $\xi^\odot=\om^\odot$.

We note further that the main theorem for the case where $E_+=E_\infty$ for the central unital unit $\om^\odot$ has been considered in \cite{Ske01}. Here we have the generalization where we assume neither that $\vt$ is a dilation of the trivial semigroup (as in \cite{Ske01}) nor of any other CP-semigroup.
\erem

\section{An open problem}

{\parskip0.5ex plus 0.5ex minus 0.5ex
\nbd{E_0}Semigroups (normal and strongly continuous in the case $\T=\R_0$) on $\sB(H)$ ($H$ a separable infinite-dimensional Hilbert space) are classified by their product system up to cocycle conjugacy; see Arveson \cite{Arv89}. This means that under the restriction on the dimension of $H$, we have a complete characterization of \nbd{E_0}semigroups up to cocycle conjugacy by product systems. (Also, under suitable technical conditions every product system of Hilbert spaces comes from an \nbd{E_0}semigroup; see Arveson \cite{Arv89a}. The corresponding question for general product systems of Hilbert modules is completely open, so far, in the continuous case $\T=\R_+$ and has been solved only recently in the discrete case $\T=\N_0$ in \cite{Ske04p}. But this is not the open problem about which we wish to speak.) The truth is that we are speaking about \nbd{E_0}semigroups only on those Hilbert spaces that are in the same isomorphism class. Under this assumption we obtain the same statement also for Hilbert modules \cite{Ske02} (without any continuity assumption): If $\vt^1$ and $\vt^2$ are strict \nbd{E_0}semigroups on $\sB^a(E^1)$ and $\sB^a(E^2)$ and $E^1\cong E^2$, then $\vt^1$ and $\vt^2$ are cocycle conjugate, if and only if their associated product systems ${E^1}^\odot$ and ${E^2}^\odot$ are \hl{isomorphic}, that is, if there exists a family $w^\odot=\bfam{w_t}_{t\in\T}$ of unitaries $w_t\in\sB^{a,bil}(E^1_t,E^2_t)$ such that $w_{s+t}=u^2_{s+t}(w_s\odot w_t)u^{1*}_{s+t}$ and $w_0=\id_\cB$. (In \cite{Ske04p} we have relaxed the condition $E^1\cong E^2$ to $\sB^a(E^1)\cong\sB^a(E^2)$. Under this condition we have cocycle conjugacy (in an obvious sense), if and only if the product systems are Morita equivalent. The notion of Morita equivalence of correspondences is borrowed from Muhly and Solel \cite{MuSo00}.)

In full generality, we do not know what we can say about the relation among $E^1$ and $E^2$ given the information that there exist \nbd{E_0}semigroups $\vt^1$ and $\vt^2$ on them that have isomorphic (or Morita equivalent) product systems. We know by explicit examples that neither of the conclusions $E^1\cong E^2$ or $\sB^a(E^1)\cong\sB^a(E^2)$ needs to be true. But what can we say, if $E^1$ and $E^2$ are inductive limits in the sense of Section \ref{ILsec} of the same product system $E^\odot$ but with respect to possibly different units ${\xi^1}^\odot$ and ${\xi^2}^\odot$? To state the problem we wish to pose clearly: Are the inductive limits over a product sytem with respect to two different units always isomorphic Hilbert modules or not?

Certainly the inductive limits will be isomorphic, if there exists an automorphism of $E^\odot$ that sends ${\xi^1}^\odot$ to ${\xi^2}^\odot$. In this case, necessarily the CP-semigroups generated by ${\xi^1}^\odot$ and ${\xi^2}^\odot$ coincide. But even in the case of Arvesons product systems of Hilbert spaces it is an open problem, whether for every pair of (normalized) units there is an automorphism of the product system that sends one unit to the other, that is, whether the automorphisms of a product system act \hl{transitively} on the set of units. A positive answer is known only for type I systems of Hilbert spaces, that is, for symmetric Fock spaces. For time ordered Fock modules we have the result provided the units generate the same CP-semigroup. (For more we cannot ask, so the statement is analogue to that for Hilbert spaces.) But while we know that the inductive limit over time-ordered Fock modules for a central unital unit (that plays the role of the vaccum) is a time-ordered Fock module (indepent of the choice of that unit), we do not know whether the same is true for a unit that generates a nontrivial CP-semigroup.

This discussion includes the representation space of the minimal weak dilation of an arbitrary nontrivial uniformly continuous unital CP-semigroup. The fact that dilations of such CP-semigroups may be obtained with help of quantum stochastic calculi lets us suspect that also the minimal weak dilation lives on a Fock module. Positive answers exist only in the case $\cB=\sB(H)$. The situation we met in the proof of the supplement when we start with a primary dilation and a spatial product system, so that there are two units arround, $\xi^\odot$ and $\om^\odot$, one of which generates a nontrivial CP-semigroup and the other unit generates the trivial one. We would be happy if we could show that the two inductive limits are isomorphic and, therefore, the minimal weak dilation a cocycle perturbation of a dilation of the trivial CP-semigroup. (This is exactly what quantum stochastic calculus usually does: Constructing a cocycle that transforms a dilation of the trivial CP-semigroup into a dilation of a nontrivial one.) We suspect that this might not be possible in general. But we have the feeling that chances might improve, when we try cocylce perturbations of $\alpha$ instead of $\vt$.
}

\setlength{\baselineskip}{2.5ex}



\newcommand{\Swap}[2]{#2#1}\newcommand{\Sort}[1]{}
\providecommand{\bysame}{\leavevmode\hbox to3em{\hrulefill}\thinspace}
\providecommand{\MR}{\relax\ifhmode\unskip\space\fi MR }
\providecommand{\MRhref}[2]{%
  \href{http://www.ams.org/mathscinet-getitem?mr=#1}{#2}
}
\providecommand{\href}[2]{#2}

\noindent
{\small\itshape Address: Dipartimento S.E.G.e S., Università degli Studi del Molise, Via de Sanctis, 86100 Campobasso, Italy. E-mail: \tt{skeide@math.tu-cottbus.de}.\\Homepage: \tt{http://www.math.tu-cottbus.de/INSTITUT/lswas/\_skeide.html}}


\end{document}